\documentclass[preprint,12pt]{article}
\usepackage{amsmath, amsthm, amssymb, amsrefs, latexsym, amsfonts, graphics, ragged2e, hyperref, bbm, 
 multicol, enumitem, xfrac}
\DeclareMathOperator*{\argmin}{\arg\!\min}
\newtheorem{theorem}{Theorem}
\newtheorem{example}[theorem]{Example}

\title{On Burdet and Johnson's Algorithm for Integer Programming}
\author{Babak Moazzez\thanks{Department of Computer Science, University of California, Davis, CA, USA} \and Kevin Cheung\thanks{School of Mathematics and Statistics, Carleton University, Ottawa, ON, Canada}}
\begin{document}
\maketitle
\begin{abstract}

 In this paper, some deficiencies of a method proposed by Burdet and Johnson in 1977 for solving integer programming problems are discussed. Examples where the algorithm fails to solve the IP and ways to fix these errors are given.
 
 \
 
\end{abstract}

\section{Introduction}\label{} 
 In 1977, Burdet and Johnson \cite{BJIP} proposed an algorithms to solve integer programs. The algorithm uses idea of lifting a subadditive dual feasible function. This function is lifted until at least one duality or subadditivity constraint gets violated. At this time, the function gets fixed for those points and lifting continues until the algorithm finds the optimal solution to IP. Although this method works for a large family of problems, it has some deficiencies and flaws which make it fail to work in general.

 Consider the Pure Integer Programming Problem (IP) and its subadditive dual:
 \begin{minipage}[t]{0.2\columnwidth}
 \begin{eqnarray}
 \begin{array}{lll}
 \label{IP}
 \mbox{min} & cx \\
 \mbox{s.t.}& Ax=b\\
  & x\geq 0, \mbox{integer}\\
 \end{array}
 \end{eqnarray}
 \end{minipage}
 \begin{minipage}[t]{0.7\columnwidth}
 \begin{eqnarray}
 \begin{array}{lll}
 \label{dual}
 \mbox{max} & F(b) \\
 \mbox{s.t.}& F(a_j)\leq c_j, \forall j  \\
  & F\in \Gamma^m , F(0)=0\\
 \end{array}
 \end{eqnarray}
 \end{minipage}
 
 \vspace{0.5cm}
 
 \noindent where $A\in\mathbb{Q}^{m\times n}$ and $b\in\mathbb{Q}^m$ and $\Gamma^m=\{F:\mathbb{R}^m\rightarrow \mathbb{R} | F \mbox{ subadditive }\}$.
 
 Gomory introduced in \cite{gomory} the group problem (corner relaxation) of an IP as:
 \begin{eqnarray}\mbox{ max }\{z^*+\bar{c_N}x_N: A_B^{-1}A_Nx_N\equiv A_B^{-1}b \ (mod\ 1), x_N\geq 0\}\label{corner}\end{eqnarray}
 where $z^*=c_BA_B^{-1}b$ and $\bar{c_N}=c_N-c_BA_B^{-1}A_N$ for basis $B$ and set of non-basic variables $N$. The non-negativity constraints on the basic variables are relaxed. Now if we add these constraints back into the problem in terms of non-basic variables, the problem will be equivalent to IP. It will be in the following form:
 \begin{eqnarray}
 \begin{array}{lll}
 \label{IPGR}
 \mbox{min} & \sum_{j\in N} \bar{c}_jx_j \\
 \mbox{s.t.}& Gx\equiv g_0 (mod\ 1)\\
 & Hx\geq h_0\\
 & x_j\geq 0, \mbox{ integer for all } j\in N.
 \end{array}
 \end{eqnarray}

 \section{The Algorithm}
  Let $S_I$ and $S_L$ denote the set of feasible points to IP (1) and its LP relaxation respectively.  If $\pi$ is a subadditive function i.e.
 $\pi(x)+\pi(y)\geq\pi(x+y) \mbox{ for all } x,y$, 
 then the inequality 
 $\sum_{j\in N} \pi_jx_j\geq \pi_0$
 is valid for (\ref{IPGR}) where
 $\pi_j=\pi(\delta^j)$, $\delta^j$ is the $j$-th unit vector and
 $\pi_0\leq \min \{\pi(x)|x\in S_I\}.$ 
 Based on the idea from \cite{BJGP} and \cite{johnson}, $\pi(x)$ is defined from a subadditive function $\Delta$ on $\mathbb{R}^N_+$ with a finite generator set $E$ by
 \[\pi(x)=\min_{y\in I(x)}\{\bar{c}y+\Delta(x-y)\}\]
 where $I(x)= E\cap S(x)$ and $S(x)=\{y\in\mathbb{Z}^N_+ | y\leq x\}$. $E\subset\mathbb{Z}_+^N$ is called the \textit{Generator Set}. Initially $E=\{0\}$. However, it will be expanded sequentially in the algorithm. Given the generator set $E$, the \textit{Candidate Set} $C$ is defined to be 
 \[C=\{x\in\mathbb{Z}_+^N: x\notin E, S(x)\backslash \{x\}\subseteq E\}.\]
 If $\pi(y_1+y_2)\leq \bar{c}y_1+\bar{c}y_2$, for all $y_1,y_2\in E$ and $y_1+y_2\in C$,
 then $\pi$ is subadditive.

 Assume that $ G $ and $ H $ have $ m_1 $ and $ m_2 $ rows respectively. For $ 2m_1+2m_2 $ real numbers
 \[\gamma^+_1,...,\gamma^+_{m_1},\gamma^-_1,...,\gamma^-_{m_1},\alpha^+_1,...,\alpha^+_{m_2},\alpha^-_1,...,\alpha^-_{m_2}\]
 Burdet and Johnson defined the generalized diamond gauge\footnote{A function which is non-negative, convex and positively homogeneous is called a gauge \cite{rock}. A generalized gauge doesn't have the non-negativity constraint.} function  $ D $ on $ \mathbb{R}^n $ to $ \mathbb{R} $ as	
 \[D(x)=\max_{\alpha,\gamma}\{\gamma Gx + \alpha Hx\}\]
 where the minimum is taken over $ 2m_1+2m_2 $ possible values: $ \gamma_i=\gamma^+_i $ or $ \gamma^-_i $ and $ \alpha_i=\alpha^+_i $ or $ \alpha^-_i $. Burdet and Johnson defined the gauge function $ \Delta_0 $ as: 
 \[\Delta_0(x)=\min_z\{D(z)| z\geq 0 , z\mbox{ integer } , Gz\equiv Gx \mbox{ (mod 1) } , Hz\geq Hx\}.\]
There are other functions used in \cite{BJIP}, but $\Delta_0$ gives the best lower bound among them. Burdet and Johnson's algorithm may be written in the following form:
  
  \vspace{0.5cm}
  
  \noindent \textbf{Initialization:} Assume that $\min_{x\in S}\{\pi(x)\}>0$ for some $S\supseteq S_I$. Scale $ \pi $ such that this minimum is equal to one. Let $ \alpha_0=\min_{j=1,...,n}\{\frac{\bar{c}_j}{\Delta(\delta^j)}|\Delta(\delta^j)>0\} $.  The initial subadditive function is $ \pi(x)=\alpha_0\Delta(x) $ and the initial bound is $ \alpha_0 $, $ E=\{0\} $ and $ C=\{\delta^i: i=1,...,n\};$
  
  \noindent \textbf{while} $ cx\neq\pi_0 $ for some  $ x $ in $ C $ and feasible to IP (\ref{IP}):
 \begin{enumerate}[itemsep=0mm]
 \item Calculate $ x^*=\argmin_{j=1,...,n} \{\frac{\bar{c}x}{\Delta(x)}|\Delta(x)>0,x\in C\} $.
 \item Let $ \alpha_0=\tfrac{\bar{c}x^*}{\Delta(x^*)} $. 
 \item Move $ x^* $ from $ C $ to $ E $ and update $ C $: adjoin to $ C $ all the points $ x>x^* $ and $ x\in \mathbb{Z}^n_+ $ with the property that for every $ y<x $ and $ y\in \mathbb{Z}^n_+ $, we have $ y\in E $.
 \item Evaluate $ \pi_0=\min_{y\in E\cap S(x)}\{\bar{c}y+\min_{x\in S} \alpha_0\Delta(x-y)\} $.
 \item (optional) Solve the following LP and update $ \pi_0 $ accordingly:
 \begin{eqnarray*}
 \begin{array}{lll}
 \mbox{max} & \pi_0 \\
 \mbox{s.t.}& \pi_0\leq \bar{c}y+\Delta(x-y)& y\in E, x\in S\\
 & \Delta(x)\leq \bar{c}x& x\in C.
 \end{array}
 \end{eqnarray*}
 \end{enumerate}

 \section{Deficiencies and ways to fix errors}
 Burdet and Johnson claimed that one could use the algorithm without using the parameter optimization. It is stated that enumeration must be done to some extent in order to proceed (which is correct), but completion of parameter adjustment is not necessary.  In the following example, we show that either with or without using parameter adjustment the algorithm fails to solve the IP. However, without parameter optimization the algorithm may fail in general.
 \begin{example}
 Consider the following integer program and its equivalent form:
 \begin{minipage}{0.5\columnwidth}
 \begin{eqnarray*}
 \begin{array}{lll}
 \mbox{min} & x_1+\ 3x_2 \\
 \mbox{s.t.}&\hspace{.9cm} -\tfrac{1}{2}x_2+x_3\hspace{0.8cm} = \tfrac{1}{2}\\
 &-x_1+x_2\hspace{0.8cm}+x_4 =-1\\
 & x_1,x_2,x_3,x_4\geq 0 \mbox{ integer.}
 \end{array}
 \end{eqnarray*}
 \end{minipage}
\begin{minipage}{0.5\columnwidth}
\begin{eqnarray*}
 \begin{array}{lll}
 \mbox{min} & x_1+\ 3x_2 \\
 \mbox{s.t.}&\hspace{.6cm} \tfrac{1}{2}x_2\equiv \tfrac{1}{2}\\
 &x_1-2x_2\geq 1\\
 & x_1,x_2\geq 0 \mbox{ integer.}
 \end{array}
 \end{eqnarray*}
\end{minipage}

\vspace{0.5cm}

The optimal solution is $ x^*=(2,1,1,0) $ with $ z^*=5 $. 
 Note that one constraint $(\tfrac{1}{2}x_2\geq -\tfrac{1}{2})$ has been removed since it is redundant. With our previous notation, we have:
 \[G=\left [\begin{array}{ll}0 & \tfrac{1}{2}\end{array}\right ],H=\left [\begin{array}{ll}1 & -1\end{array}\right ], g_0=\tfrac{1}{2}\mbox{ and }h_0=1.\]
 Choose $\Delta_0$ with $ S=S_L$. Scale $ \gamma $ and $ \alpha $ so that $ \min_{x\in S} \{\pi(x)\}=1 $. Note that $x=(1,0)$ will minimize $\Delta_0$ on $S$ since any point $x$ that we choose from $S$, will have $x_1-x_2\geq 1$. So we have:
 \begin{eqnarray*}1&=& \min_z\{\max \left\{\begin{array}{c} \tfrac{1}{2}\gamma^+ z_2\\ -\tfrac{1}{2}\gamma^- z_2 \end{array}\right\} + \max \left\{\begin{array}{c} \alpha^+ z_1-\alpha^+z_2\\ -\alpha^- z_1+\alpha^-z_2 \end{array}\right\}:\\ & &z\geq 0 ,\mbox{integer} , \tfrac{1}{2}z_2\equiv 0 , z_1-z_2\geq 1\}.\end{eqnarray*}
 
 \noindent $z_2$ can take values $2k$ for $k=0,1,2,...$ and $z_1=1+2k$. Obviously $k=0$, will give the optimal value and this gives us $\alpha^+=\alpha^-=1$. $\gamma$ could be any value. We choose $\gamma^+=\gamma^-=1$. 
  
  Note that for $ k\geq 0 $, $ \alpha^+ = \alpha^- =\alpha $ and $ \gamma^+ = \gamma^- = \gamma $ for some $ \alpha $ and $ \gamma $, we have 
\[ \Delta_0(k,0) =\min_z\{ \max \left\{\begin{array}{c} \alpha z_1-\alpha z_2\\ -\alpha z_1+\alpha z_2 \end{array}\right\}:z\geq 0 ,\mbox{integer} ,z_1-z_2\geq k\}=k\alpha. \]
   Now $ E=\{(0,0)\} $ and $ C=\{(1,0),(0,1)\} $. 
 We have $ \Delta_0(1,0)=1 $ and $ \Delta_0(0,1)=\tfrac{1}{2}$, so
 \[\alpha_0=\min_{j=1,...,n}\{\tfrac{\bar{c}_j}{\Delta_0(\delta^j)}|\Delta_0(\delta^j)>0\}=\min \{\tfrac{1}{\Delta_0(1,0)},\tfrac{3}{\Delta_0(0,1)}\}=\min \{\tfrac{1}{1},\tfrac{3}{\tfrac{1}{2}}\}=1,\]
 so initial $ \pi(x) $ is $\Delta_0(x)$. 
 \[x^*=\argmin_{j=1,...,n} \{\tfrac{\bar{c}x}{\Delta_0(x)}|\Delta_0(x)>0,x\in C\}=\argmin \{\tfrac{1}{\Delta_0(1,0)},\tfrac{3}{\Delta_0(0,1)}\} =(1,0).\]
 Now the new sets are $ E=\{(0,0),(1,0)\} $ and $ C=\{(2,0),(0,1)\} $.
 \[\pi_0=\min_{y\in E\cap S(x)}\{\bar{c}y+\min_{x\in S} \alpha_0\Delta_0(x-y)\}=1. \]
 
 At this point if we continue without parameter optimization, $ \pi_0 $ will always remain at $ 1 $ and the points $ (k,0) $ for $ k\geq 2 $ will enter $ C $ and then $ E $ one by one. As result, $ (2,1) $ which corresponds to the optimal solution will never enter $ C $ and hence the algorithm will never find the optimal solution. This can be seen by noting that $ \Delta_0(k,0)=k $  and $ \bar{c}\cdot(k,0)=k $ for $ k\geq 2 $, while $ \Delta_0(0,1)=\tfrac{1}{2} $ and $ \bar{c}\cdot (0,1)=3 $.
 
 Otherwise if we use parameter adjustment, the LP will have the form
 \begin{eqnarray*}
 \begin{array}{lll}
 \mbox{max} &\pi_0 \\
 \mbox{s.t.}&\pi_0 \leq \Delta_0(x) & x\in S_L\\
 &\pi_0 \leq 1+ \Delta_0(x-(1,0)) & x\in S_L \\
 & \Delta_0(0,1)\leq 3, \Delta_0(2,0) \leq 2
 \end{array}
 \end{eqnarray*}
 
 \noindent Note that the first constraint can be written as $ \pi_0\leq \alpha^+$.
 Also we have $ \Delta_0(0,1)=\max \{ \tfrac{1}{2}\gamma^+ , -\tfrac{1}{2}\gamma^- \}$ and $\Delta_0(2,0)= 2\alpha^+$. So the parameter adjustment LP becomes equivalent to 
  \begin{eqnarray*}
   \begin{array}{lll}
   \mbox{max} &\pi_0 \\
   \mbox{s.t.}&\pi_0 \leq\min\{1, \alpha^+\}\\
   & \max \{ \tfrac{1}{2}\gamma^+, -\tfrac{1}{2}\gamma^- \} \leq 3\\
   & 2\alpha^+ \leq 2.
   \end{array}
   \end{eqnarray*}
  This LP gives $ \alpha^+=\alpha^-=1, \gamma^+=\gamma^-=6 $ and $ \pi_0=1$.  With these new parameters, we will have:
 $\Delta_0(0,1)=3,  \Delta_0(2,0)=2$ and $x^*=(0,1)$.
  New $ E=\{(0,0),(1,0),(0,1)\} $ and $ C=\{(2,0),(1,1),(0,2)\} $. $ \alpha_0=1 $ and $\pi_0=1$.
 At this point if one continues without parameter optimization, since $ \Delta_0(0,2)=0 $, the points $ (k,0) $ for $ k\geq 2 $ will enter $ C $ and then $ E $ one by one and $ (2,1) $  will never enter $ C $. So again we continue with parameter optimization. 
  The parameter adjustment LP becomes equivalent to 
   \begin{eqnarray*}
    \begin{array}{lll}
    \mbox{max} &\pi_0 \\
    \mbox{s.t.}&\pi_0 \leq \min\{\alpha^+,1\}\\
    & \max \{ \tfrac{1}{2}\gamma^+, -\tfrac{1}{2}\gamma^- \} \leq 4\\
    & 2\alpha^+ \leq 2.
    \end{array}
    \end{eqnarray*}
 This LP gives $ \alpha^+=\alpha^-=1, \gamma^+=\gamma^-=8 $ and $ \pi_0=0$.  With these new parameters, we will have:
 \[\Delta_0(1,1)=4 , \Delta_0(0,2)=0, \Delta_0(2,0)=2 \mbox{ and } x^*=(1,1).\]
 New $ E=\{(0,0),(1,0),(0,1),(1,1)\} $ and $ C=\{(2,0),(0,2)\} $. $ \alpha_0=1 $ and $\pi_0=1$.
 First Note that for any $ \alpha $  and $ \gamma, \Delta_0(0,2)=0 $.
 From this iteration on, since $ (2,0) $ is in $ C $ and $ \Delta_0(0,2)=0 $, again the points $ (k,0) $ for $ k\geq 3 $ will enter $ C $ and then $ E $ one by one and $ (2,1) $  will never enter $ C $ and this algorithm will never terminate. Since we have in the constraints of parameter adjustment LP that $ \Delta_0(k,0)\leq k $ which gives $ k\alpha^+\leq k $  and $ \alpha^+\leq 1 $, $ \alpha^+ $ will never change. Consequently, $ \pi_0 $ will never increase to more than 1 with these settings because of the first constraint in parameter adjustment LP namely  $ \pi_0\leq \alpha^+ $ which will remain in the LP in all iterations since $ (0,0) \in E $ in all iterations. In iteration $ k\geq 3 $, the LP will have the following form:
    \begin{eqnarray*}
     \begin{array}{lll}
     \mbox{max} &\pi_0 \\
     \mbox{s.t.}&\pi_0 \leq \min\{\alpha^+,1\}\\
     & \max \{ \tfrac{1}{2}\gamma^+, -\tfrac{1}{2}\gamma^- \} \leq k+1\\
     & (k+1)\alpha^+ \leq k+1.
     \end{array}
     \end{eqnarray*}
 \end{example}
 
 This problem can be fixed by imposing an upper bound for each variable. In this case, eventually $\pi_0$ and $\alpha_0$ will increase. For instance in the example above,
 if we have $x_1\leq K$, then $(0,2)$ will enter $C$ if we get to a point that we have
 $E=\{(0,0),(1,0),...,(K,0)\}$ and  $C=\{(K+1,0),(0,2)\}.$ 
 This means that entering  $(K+1,0)$ to $E$ will not help since it does not belong to $X_I$.
 See \cite{schrijver} for imposing bounds on variables.
 
 The number of steps for the enumeration part of the algorithm is proportional to the coordinates of the optimal solution. For example if $ x^* $ is optimal solution for some IP with $ x^*_i=1000 $ for some $ i $, it will take 1000 iterations to find the optimal solution. The reason is that $ x^* $ must enter $ C $ at some iteration and $ \pi_0 $ has to increase to $ \bar{c}x^* $. Even if we move multiple points from $ C $ to $ E $ at each iteration, the $ x_i $ for the points in $ C $ will increase only by one unit at most.  However, this algorithm may work better for binary problems. 
  
 The following example shows that during the algorithm, strong duality may not hold all the time. 
 \begin{example}\label{BJfail}
 Consider the following integer program and its equivalent form:
 \begin{minipage}{0.4\columnwidth}
 \begin{eqnarray*}
 \begin{array}{lll}
 \mbox{min} & x_1+\ c_2x_2 \\
 \mbox{s.t.}&x_1+\tfrac{1}{2}x_2+x_3 = \tfrac{1}{2}\\
 & x_1,x_2,x_3\geq 0 \mbox{ integer.}
 \end{array}
 \end{eqnarray*}
 \end{minipage}
\begin{minipage}{0.6\columnwidth}
\begin{eqnarray*}
 \begin{array}{lll}
 \mbox{min} & x_1+\ c_2x_2 \\
 \mbox{s.t.}&\hspace{1.1cm} \tfrac{1}{2}x_2\equiv \tfrac{1}{2}\\
 & -x_1-\tfrac{1}{2}x_2\geq -\tfrac{1}{2} \\
 & x_1,x_2\geq 0 \mbox{ integer.}
 \end{array}
 \end{eqnarray*}
\end{minipage}

\vspace{0.5cm}

 Choose any subadditive function with $\Delta(0)=0$ and $\Delta(\delta^1)>0$. Here, $S_I=\{\delta^2\}$ and $z_{IP}=c_2$.
 For any $\alpha\geq 0$ we have
 \[\pi_0=\min_{x\in S_I}\pi(x)=\min_{y\in E,y\leq \delta^2}\{\bar{c}y+\alpha\Delta(\delta^2-y)\}=\min\{\alpha\Delta(\delta^2),c_2\}.\]
 Since $\alpha\Delta(\delta^1)\leq c_1$, we get that $\alpha\leq \frac{1}{\Delta(\delta^1)}$. If $c_2>\frac{\Delta(\delta^2)}{\Delta(\delta^1)}$, then 
 $\pi_0=\frac{\Delta(\delta^2)}{\Delta(\delta^1)}<c_2=z_{IP}$. So strong duality does not hold.
 \end{example}
 This problem can be fixed by performing some preprocessing first. 
 If $ \delta_i \notin X_I $ then $ x_i $ must be eliminated from IP by letting it be zero.

\bibliographystyle{model1a-num-names}

\end{document}